# Нестационарные моменты для систем обслуживания

Nonstationary moments for queuing systems

**Головастова Э.А.**

## Аннотация

В данной работе доказана экспоненциальная сходимость нестационарного момента случайной величины, задающей виртуальное время ожидания в системе $M\backslash G\backslash 1\backslash\infty$, к моменту виртуального времени ожидания при стационарном режиме.

**Ключевые слова:** условие Крамера, скорость сходимости, математическое ожидание, виртуальное время ожидания.

## Результаты для общего случая

Рассмотрим систему массового обслуживания типа $M\backslash G\backslash 1\backslash\infty$ с коэффициентом загрузки $\rho = \lambda b_1 < 1$, при дисциплине обслуживания - FIFO. Входящий поток - пуассоновский с параметром $\lambda$, время обслуживания распределено по закону $B(x)$. Полагаем, что функция $B(x)$ обладает всеми необходимыми свойствами для справедливости последующих преобразований. В частности, у $B(x)$ существуют все моменты. Обозначим $b_k = \int_{-\infty}^{\infty} x^k dB(x)$ ее $k$-ый момент.

Пусть $\beta(s) = \int_0^{\infty} e^{-sx} dB(x)$ - преобразование Лапласа-Стилтьеса функции $B(x)$. Также положим, что $\overline{B(x)} = 1 - B(x) \leq Ce^{-\gamma x}$, $\gamma > 0$, то есть распределение времени обслуживания имеет экспоненциальный хвост. Последнее эквивалентно следующему определению:

**Определение.** Если $\exists \delta_0 > 0$ такое, что

$$\int_0^\infty e^{\delta_0 x} dB(x) = \beta(-\delta_0) < \infty. \quad (1)$$

То говорят, что функция $B(x)$ удовлетворяет условию Крамера.

**Утверждение.** Период занятости в рассматриваемой системе удовлетворяет условию Крамера.

*Доказательство.* Пусть случайная величина $\tau$ - период занятости, $\pi(s) = Ee^{-s\tau}$.

Известно следующее уравнение [1]:

$$\pi(s) = \beta\big(s + \lambda(1 - \pi(s))\big).$$

Тогда:

$$-E\tau = \pi'(0) = \beta'\big(\lambda(1 - \pi(0))\big)(1 - \lambda\pi'(0)); \quad \pi(0) = 1$$

$$\pi'(0) = \frac{\beta'(0)}{\lambda\beta'(0) + 1} \Rightarrow E\tau = \frac{b_1}{1 - \lambda b_1}.$$

Заметим, что $\int_0^\infty e^{sx} dB(x) = \pi(-s)$ монотонно возрастает.

Тогда покажем, что $\exists \delta_1 > 0$ такое, что $\pi(-s) < \infty \quad \forall s < \delta_1$.

Пусть $g(s) = \pi(-s) = \beta\big(-s + \lambda(1 - g(s))\big)$, $g'(0) = E\tau$. Тогда:

$$g(s) = \beta\big(-s + \lambda(1 - 1 - g'(0)s) + o(s^2)\big) = \beta\big(-\frac{s}{1 - \lambda b_1} + o(s^2)\big).$$

Для $\beta(s)$ выполнены условия Крамера, т.е. $\exists \delta_0 > 0$ такое, что $\beta(-\delta_0) < \infty$. Тогда при $s$: $\frac{s}{1 - \lambda b_1} < \delta_0 \Leftrightarrow s < \delta_0(1 - \lambda b_1) = \delta_1$ выполнено:

$$\pi(-s) = g(s) = \beta\big(-\frac{s}{1 - \lambda b_1} + o(s^2)\big) < \infty \qquad \square$$

Из утверждения следует, что период занятости имеет экспоненциальный хвост в рассматриваемой системе.

Обозначим через $W(t)$ - виртуальное время ожидания, то есть то время, которое пришлось бы ожидать начала обслуживания требованию, поступившему в момент времени $t$. Это регенерирующий процесс, моменты регенерации которого есть моменты попадания $W(t)$ в $0$, или, по-другому, моменты начала периодов занятости. Считаем $W(0) = 0$ почти наверное.

Пусть $\phi(t) = EW(t)$. Период регенерации процесса $W(t)$ - случайная величина $\zeta_i$ - есть сумма периода занятости и периода свободного состояния системы. Эти два периода - независимые случайные величины, функции распределения которых удовлетворяют условию Крамера. Поэтому функция распределения $\zeta_i$ также удовлетворяет условию Крамера, то есть имеет экспоненциальный хвост.

Пусть $H(t)$ - функция восстановления, соответствующая последовательности $\{\zeta_i\}$. Тогда из вышесказанного следует, что процесс $W(t)$ удовлетворяет условиям результата из [2] (см. Теорема 3. Приложение1.). Использовав его, имеем:

$$\exists \varepsilon > 0: \quad H(t) = \frac{t}{\mu} + \frac{\mu_2}{2\mu^2} + R(t); \quad |R(t)| \leq c e^{-\varepsilon t}, \tag{2}$$

$$\mu = E\zeta_i = \frac{1}{\lambda} + \frac{b_1}{1 - \lambda b_1} = \frac{1}{\lambda(1 - \lambda b_1)}, \quad \mu_2 = E(\zeta_i)^2$$

Также известен следующий результат:

$$h(t) = H'(t) = \frac{1}{\mu} + r(t); \quad |r(t)| \leq C_1 e^{-\varepsilon_1 t}$$

Далее, так как $W(t)$ - марковский процесс, то можно написать следующее:

$$P(W(t) > x) = P(W(t) > x, \zeta_1 > t) + \int_0^t P(W(t-y) > x) dP(\zeta_1 < y)$$

$$\phi(t) = \int_0^\infty P(W(t) > x) dx = \int_0^\infty P(W(t) > x, \zeta_1 > t) dx + \int_0^t \left( \int_0^\infty P(W(t-y) > x) dx \right) dP(\zeta_1 < y)$$

$$\phi(t) = q(t) + \int_0^t \phi(t-y) dF(y), \tag{3}$$

где

$$F(y) = P(\zeta_1 < y), \quad q(t) = E\big(W(t)\mathbb{I}(\zeta_1 > t)\big).$$

Последнее есть уравнение восстановления. Его решение:

$$\phi(t) = \int_0^t q(t-y) dH(y) = \int_0^t q(t-y) h(y) dy \tag{4}$$

Также известен следующий факт [3]:

$$\lim_{t \to \infty} \phi(t) = \frac{1}{\mu} \int_0^\infty q(y) dy = \frac{\lambda b_2}{2b_1(1-\rho)} = \phi$$

**Теорема.** $\phi(t)$ сходится к $\phi$ экспоненциально быстро.

*Доказательство.* Наша цель - оценить $\phi(t) - \phi$. Из (4) имеем:

$$\phi(t) - \phi = \int_0^t q(t-y)h(y)dy - \frac{1}{\mu}\left(\int_0^t + \int_t^\infty\right)q(y)dy$$

$$|\phi(t) - \phi| = \left|\int_0^t q(t-y)h(y)dy - \frac{1}{\mu}\int_0^t q(y)dy\right| + \frac{1}{\mu}\int_t^\infty q(y)dy = J_1 + J_2$$

Пусть $\gamma_t = (\zeta_1 - t)\mathbb{I}(\zeta_1 > t)$. Тогда:

$$W(t)\mathbb{I}(\zeta_1 > t) \leq (\zeta_1 - t)\mathbb{I}(\zeta_1 > t) \leq \zeta_1 \mathbb{I}(\zeta_1 > t)$$

Также вспомним, что $\exists \varepsilon_2 > 0 : \overline{F(x)} \leq e^{-\varepsilon_2 t}$.

Тогда:

$$q(t) = E\left(W(t)\mathbb{I}(\zeta_1 > t)\right) \leq -\int_t^\infty x\, d\overline{F(x)} = -x\overline{F(x)}\Big|_t^\infty + \int_t^\infty \overline{F(x)}\, dx = t\,\overline{F(t)} + \int_t^\infty \overline{F(x)}\, dx \leq C_3 e^{-\varepsilon_3 t}$$

для некоторого $\varepsilon_3$ и $\forall t$. Отсюда получается экспоненциальная оценка для $J_2$.

Далее воспользуемся выражением для $h(y)$, получим:

$$\int_0^t q(t-y)h(y)\, dy = \frac{1}{\mu}\int_0^t q(y)\, dy + \int_0^t q(t-y)r(y)\, dy.$$

Тогда, используя полученную оценку для $q(t)$ и то, что $|r(t)| \leq C_1 e^{-\varepsilon_1 t} < 1$ для $\forall t > 0$, имеем:

$$J_1(t) \leq \int_0^t q(t-y)|r(y)|\, dy = \left(\int_0^{t/2} + \int_{t/2}^t\right)q(t-y)|r(y)|\, dy \leq C_1\int_{t/2}^t q(y)dy + C_3\int_{t/2}^t |r(y)|\, dy \leq C_4 e^{-\varepsilon_4 t}$$

$\square$

## Частный случай

Проверим полученный выше результат для системы $M\setminus M\setminus 1\setminus \infty$ при дисциплине обслуживания - FIFO. Также полагаем $\rho = \frac{\lambda}{\mu} < 1$. Здесь входящий поток - пуассоновский с параметром $\lambda$, время обслуживания распределено по закону $B(x) = 1 - e^{-\mu x}, x \geq 0$; $\beta(s) = \frac{\mu}{\mu + s}$. Для этой системы имеется явное выражение для распределения виртуального времени ожидания [4]:

$$W(t, x) = P_0(t) + \sum_{k=1}^\infty P_k(t)\, P\left(\sum_{i=1}^k \xi_i < x\right),$$

где $P_n(t)$ - вероятность того, что в момент $t$ в системе $n$ требований, $\xi_i$ - независимые случайные величины, распределенные экспоненциально: $exp(\mu)$. Учтем, что $Y_k = \sum_{i=1}^{k} \xi_i$, поэтому имеет распределение $\Gamma(\frac{1}{\mu}, k)$. Тогда:

$$\phi(t) = EW(t,x) = P_0(t) + \sum_{k=1}^{\infty} P_k(t) C_k; \quad C_k = EY_k = \frac{k}{\mu}. \quad (5)$$

Положим $P_n(0) = \delta_n^0$. Для $P_n(t)$ также известно явное выражение [1]:

$$P_n(t) = e^{-(\lambda+\mu)t}\left[(\sqrt{\frac{\mu}{\lambda}})^{-n} I_n(2\sqrt{\lambda\mu t}) + (\sqrt{\frac{\mu}{\lambda}})^{-n+1} I_{n+1}(2\sqrt{\lambda\mu t}) + (1-\frac{\lambda}{\mu})(\frac{\lambda}{\mu})^n \sum_{k=n+2}^{\infty} (\sqrt{\frac{\mu}{\lambda}})^k I_k(2\sqrt{\lambda\mu t})\right]$$

Тут $I_n(z)$ - функция Бесселя первого рода.

Преобразуем выражение для $P_n(t)$, используя следующее разложение:

$$e^{(\frac{x}{2})(y+\frac{1}{y})} = \sum_{n=-\infty}^{\infty} y^n I_n(x) = \sum_{n=0}^{\infty} y^n I_n(x) + \sum_{n=1}^{\infty} y^{-n} I_n(x).$$

Взяв в нем $y = \sqrt{\frac{\mu}{\lambda}}$ и $x = \sqrt{\lambda\mu t}$ и учитывая $I_{-n}(x) = I_n(x)$, получим:

$$P_n(t) = e^{-(\lambda+\mu)t}\left[(\sqrt{\frac{\mu}{\lambda}})^{-n} I_n(t) + (\sqrt{\frac{\mu}{\lambda}})^{-n+1} I_{n+1}(t) - \rho^n(1-\rho)\sum_{p=0}^{n+1}(\sqrt{\frac{\mu}{\lambda}})^p I_p(t)\right] +$$

$$+ e^{-(\lambda+\mu)t}\rho^n(1-\rho)exp\left(\frac{2\sqrt{\lambda\mu t}}{2}(\sqrt{\frac{\mu}{\lambda}}+\sqrt{\frac{\lambda}{\mu}})\right) -$$

$$- e^{-(\lambda+\mu)t}\rho^n(1-\rho)\sum_{p=1}^{\infty}(\sqrt{\frac{\lambda}{\mu}})^p I_p(t).$$

Известна следующая асимптотика при $t \to \infty$:

$$I_n(2\sqrt{\lambda\mu t}) \approx \frac{exp(2\sqrt{\lambda\mu t})}{\sqrt{2\pi(2\sqrt{\lambda\mu t})}}.$$

Тогда:

$$P_n(t) \approx \rho^n(1-\rho) + \frac{e^{(2\sqrt{\lambda\mu}-(\lambda+\mu))t}}{\sqrt{4\pi\sqrt{\lambda\mu t}}}\left[(\sqrt{\frac{\mu}{\lambda}})^{-n}(1+\sqrt{\frac{\mu}{\lambda}}) - \rho^n(1-\rho)\left(\frac{(\sqrt{\frac{\mu}{\lambda}})^{n+2}-1}{\sqrt{\frac{\mu}{\lambda}}-1} + \frac{\sqrt{\rho}}{1-\sqrt{\rho}}\right)\right].$$

Подставим полученное асимптотическое выражение в (5):

$$\phi(t) \approx (1-\rho) + \frac{(1-\rho)}{\mu}\sum_{k=1}^{\infty} k\rho^k + \frac{exp(\cdots)}{\sqrt{4\pi\sqrt{\lambda\mu}t}} \cdot$$

$$[\rho(1+\sqrt{\frac{\mu}{\lambda}}) - \sqrt{\frac{\lambda}{\mu}}(1+\sqrt{\frac{\lambda}{\mu}}) + \frac{1+\sqrt{\frac{\mu}{\lambda}}}{\mu}\sum_{k=1}^{\infty} k(\sqrt{\rho})^k - \frac{1-\rho}{\mu}\sum_{k=1}^{\infty} k\rho^k \big(\frac{(\sqrt{\frac{\mu}{\lambda}})^{k+2}-1}{\sqrt{\frac{\mu}{\lambda}}-1} + \frac{\sqrt{\rho}}{1-\sqrt{\rho}}\big)].$$

Тогда:

$$\phi(t) \approx (1-\rho) + \frac{\rho}{\mu(1-\rho)} + \frac{e^{(2\sqrt{\lambda\mu}-(\lambda+\mu))t}}{\sqrt{4\pi\sqrt{\lambda\mu}t}} \frac{(1+3\sqrt{\rho}+4\rho+4\rho^{3/2}+3\rho^2-\rho^{5/2})}{\mu(1-\sqrt{\rho}-\rho^2+\rho^{5/2})}.$$

Заметим, что при $\rho = \frac{\lambda}{\mu} < 1$ $\lambda + \mu > 2\sqrt{\lambda\mu}$, и тогда слагаемое с экспонентой в последнем выражении стремится к нулю при $t \to \infty$.